\documentclass[11pt,letterpaper,reqno]{amsart}
\usepackage{amssymb,amsmath,amscd,color,amsfonts}

\usepackage{verbatim}
\usepackage[all]{xy}
\usepackage[colorlinks=true]{hyperref}
\usepackage{enumerate}
\usepackage{tensor}

\newtheorem{introthm}{Theorem}  
\newtheorem{introcor}{Corollary}  
\newtheorem{theorem}{Theorem}[section]

\newtheorem{corollary}[theorem]{Corollary}
\newtheorem{lemma}[theorem]{Lemma}
\newtheorem{proposition}[theorem]{Proposition}
\newtheorem{example}[theorem]{Example}
\newtheorem{remark}[theorem]{Remark}

\theoremstyle{definition}
\newtheorem{definition}[theorem]{Definition}

\def\H{\mathcal{H}}

\def\C{\mathbb{C}}
\def\R{\mathbb{R}}

\def\G{\mathcal{G}}
\def\H{\mathcal{H}}
\def\V{\mathcal{V}}
\def\N{\mathcal{N}}

\def\O{\mathcal{O}}
\def\D{\mathcal{D}}
\def\d{\mathrm{d}}
\def\pr{\mathrm{pr}}
\def\tto{\rightrightarrows}

\DeclareMathOperator{\rank}{rank} 
\DeclareMathOperator{\Ker}{ker} 
\DeclareMathOperator{\Lie}{Lie} 
\newcommand{\diffto}{\xrightarrow{\raisebox{-0.2 em}[0pt][0pt]{\smash{\ensuremath{\sim}}}}} 

\title[Analytic Linearization and Holomorphic Extensions]{Analytic Linearization and Holomorphic Extensions of Proper Groupoids}

\author{Rui Loja Fernandes}
\address{Department of Mathematics, University of Illinois at Urbana-Champaign, 1409 W. Green Street, Urbana, IL 61801 USA}
\email{ruiloja@illinois.edu}

\author{Ning Jiang}
\address{Department of Mathematics, University of Illinois at Urbana-Champaign, 1409 W. Green Street, Urbana, IL 61801 USA}
\email{njiangims@icloud.com}

\thanks{The authors were partially supported by NSF grant DMS-2303586.}
\begin{document}

\begin{abstract}
    We establish analytic linearization of s-proper analytic grou\-poids around invariant sub\-ma\-ni\-folds. We apply this result to show that any such groupoid admits a holomorphic extension.  
\end{abstract}

\maketitle

{
    \hypersetup{linkcolor=black} 
    \setcounter{tocdepth}{1}
    \tableofcontents
}

\section{Introduction}
Proper Lie groupoids are known to be linearizable around their orbits and, more generally, around saturated submanifolds. This was conjectured by Alan Weinstein \cite{weinstein2000linearization,weinstein2002linearization} and later proved in full generality in \cite{crainic2013linearization,zung2006proper} (see also \cite{pflaum2014geometry,del2018riemannian}). All these works were carried out in the smooth category. In this paper, we show that analogous linearization results hold in the real analytic category. We apply these results to study holomorphic extensions of real analytic algebroids and groupoids.

An indication that such linearization results may exist is provided by the classical Bochner linearization theorem, which states that an action of a compact Lie group can be linearized around a fixed point (see, e.g., \cite{duistermaat2012lie}). This theorem holds in both the smooth and analytic settings and can be proved via a simple averaging argument: one first constructs a metric invariant under the group action and then applies the exponential map to achieve linearization. This approach was extended to the setting of proper Lie groupoids in \cite{del2018riemannian}, where the authors introduced the notion of a 2-metric on a Lie groupoid. They proved that (i) every proper Lie groupoid admits a 2-metric, and (ii) a Lie groupoid equipped with a 2-metric can be linearized around an orbit via the exponential map of the 2-metric.

It is not hard to see that if a real analytic groupoid admits an analytic 2-metric, then the exponential map yields an analytic linearization around an orbit. In this paper, we construct an analytic 2-metric on any real analytic s-proper Lie groupoid. The construction relies on an averaging procedure adapted to the analytic setting, which in turn requires the existence of an analytic Haar density. In Appendix \ref{sec:analytic_haar_density}, we establish that any s-proper real analytic groupoid admits such a density. Hence, the exponential map provides an analytic linearization, yielding our first main result:

\begin{introthm}\label{thm:main_theorem_1}
Every s-proper real analytic groupoid admits an analytic linearization around any invariant submanifold.
\end{introthm}

Note that it is shown in \cite{Martinez20} that every proper Lie groupoid admits a compatible real analytic structure.

Theorem \ref{thm:main_theorem_1} can be applied to promote known linearization results from the smooth to the real analytic category. One example is the following real analytic version of \cite[Theorem 8.6]{FM24}:

\begin{introcor}
Let $(M,\pi)$ be a real analytic Poisson manifold and $S\subset M$ a Poisson submanifold. If $T^*_SM$ is integrable by a compact, Hausdorff, Lie groupoid whose source fibers have trivial 2nd de Rham cohomology, then $(M,\pi)$ is analytically linearizable around $S$.
\end{introcor}

Let us turn now to holomorphic extensions. By a {\bf holomorphic extension of a real analytic algebroid} $A\Rightarrow M$ we mean a holomorphic Lie algebroid $A_\C\Rightarrow M_\C$  together with a totally real algebroid embedding $\iota:A\to A_\C$ (cf.~Definition \ref{def:holomorphic:extension:algebroid}). It is not hard to see that any real analytic Lie algebroid admits a holomorphic extension, that the germ of this extension is unique, and that it satisfies a universal property. We discuss these and other results about holomorphic extensions of real analytic algebroids in Section \ref{sec:algebroids:holomorphic:extensions}. 

At the groupoid level, the notion of holomorphic extension is more subtle. In Section \ref{sec:groupoids:holomorphic:extensions}, we introduce a notion {\bf holomorphic extension of a real analytic groupoid} which extends both the notion of complexification of a manifold and of holomorphic extension of a Lie group. Note that a real analytic groupoid may not have a holomorphic extension. This is already the case for Lie groups. However, by a classical result of Chevalley, every compact Lie group admits a holomorphic extension, and our second main theorem extends this fact.

\begin{introthm}\label{thm:main_theorem_2}
Every real analytic s-proper groupoid admits a holomorphic extension.
\end{introthm}

Our proof of this result, given in Section~\ref{sec:s:proper}, uses Theorem~\ref{thm:main_theorem_1} to obtain local holomorphic extensions around orbits, which are then glued into a global holomorphic extension.

In general, one may ask:
\begin{itemize}
\item Given an \emph{integrable} real analytic algebroid, does it admit a holomorphic extension that is also integrable?
\end{itemize}

This question is relevant, for example, to the problem of integrating complex Lie algebroids -- see the speculations about this issue by Weinstein in \cite{weinstein2007integration}. Note that the integrability criterion from \cite{crainic2003integrability} is difficult to apply in this context since it requires computing the monodromy groups of orbits. However, it is easy to construct an example of an integrable analytic Lie algebroid $A\Rightarrow M$ whose complexification $A_\C\Rightarrow M_\C$ has orbits in $M_\C-M$ with closures intersecting $M$.

Theorem~\ref{thm:main_theorem_2} provides the following partial answer to this question.

\begin{introcor}
If a real analytic algebroid $A$ admits an s-proper integration then $A$ admits a holomorphic extension which is integrable.
\end{introcor}

There are also interesting consequences of these results to Poisson geometry. For example, we have the following application to Poisson manifolds of s-proper type, i.e., which integrate to some s-proper symplectic groupoid (see \cite{CFM1,CFM2,CFM3} for a detailed discussion of this class of Poisson manifolds).

\begin{introcor}
Every real analytic Poisson manifold of s-proper type admits an extension to an integrable holomorphic Poisson manifold.
\end{introcor}

The techniques used in this paper rely on s-properness and do not extend to the general proper case. It remains unclear to us to what extent the theory developed here can be extended to the proper case. Actions of non-compact Lie groups (and non-proper Lie groupoids) can also be linearized under certain conditions; see \cite{Bischoff24,Bischoff25}. It is therefore possible that the theory can be developed also in the non-proper case under appropriate assumptions.
\smallskip 

{\bf Acknowledgments.} We would like to thank Ana Balibanu, Matias del Hoyo,  Brent Pym, Florian Zeiser, and Maarten Mol for several valuable discussions that helped shape the ideas in this paper. We are indebted to David Martínez Torres for his many insightful comments and suggestions on an earlier draft of this paper. We also thank the anonymous referees for their valuable remarks, which helped improve the presentation of the paper.

\section{Analytic 2-metrics}

{\bf Convention.} Throughout this paper we assume that Lie groupoids are Hausdorff and source connected. Also, by an analytic algebroid/groupoid we mean a real analytic one.
\smallskip 

The main aim of this section is to prove the following result:

\begin{theorem}
Every analytic s-proper groupoid admits an analytic 2-metric.
\end{theorem}

The proof to be given in the next subsections is an adaptation to the analytic setting of the proof in the smooth setting given in \cite{del2018riemannian}.

\subsection{Analytic Riemannian submersions}
The classical results of Grauert \cite{Grauert58} and Morrey \cite{Morrey58}, show that any analytic manifold can be embedded analytically in euclidean space. Hence, every analytic manifold admits analytic Riemannian metrics.

Let $(E, \eta^E)$ and $(B, \eta^B)$ be Riemannian manifolds. Recall that a map $p: E\rightarrow B$ is called \emph{Riemannian submersion} if for any $e\in E$, the differential of $p$ at $e$ restricts to an isometry between the horizontal space at $e$ and the tangent space of $B$ at $b=p(e)$:
\[\d_{e}p: (\Ker \d p)_e^{\perp} \diffto T_{b}B\]

Given any submersion $p: E\rightarrow B$ we say that a Riemannian metric $\eta^E$ is {\bf $p$-transverse} if for any two points $e,e'\in E$ with $p(e)=p(e')$, the composition
\[ (\Ker \d p)_e^{\perp} \rightarrow T_{b}B \leftarrow (\Ker \d p)_{e'}^{\perp}  \]
is an isometry. When $\eta^E$ is $p$-transverse, it induces the push-forward metric $\eta^B = p_{*}\eta^E$ on $B$, which makes $p: E\rightarrow B$ a Riemannian submersion. If we further assume that $\eta^E$ is analytic and $p$ is analytic, then the push-forward metric $\eta^B=p_* \eta^E$ is also analytic.

\begin{lemma}
\label{lem:analytic:p:transv}
    Let $p: E\rightarrow$ B be an analytic submersion. There exists a $p$-transverse analytic metric $\eta_E$ on $E$.
\end{lemma}
\begin{proof}
Choose an auxiliary analytic metric $\bar{\eta}_E$ on $E$. The orthogonal projection $TE \rightarrow (\Ker(\d p))^\perp$ gives an analytic splitting $\sigma$ of the exact sequence over $E$
    \[ \xymatrix{
    0\ar[r] & \Ker(\d p) \ar[r] & TE \ar[r]^{p_*} & p^*(TB) \ar[r] \ar@/^/@{-->}[l]^{\sigma} & 0.}
\]
This induces an analytic isomorphism 
\[TE \cong p^*TB \oplus \Ker(\d p),\]
so choosing an analytic metric $\eta_B$, we can define an analytic $p$-transverse metric by
\[ \eta_E := \bar{\eta}_E\mid_{\Ker(\d p)} + p^*\eta_B. \]
\end{proof}

\begin{lemma}[\cite{del2018riemannian}]
Let $q:\tilde{E}\rightarrow E$ and $p:E\rightarrow M$ be surjective submersions. Let $\eta$ be a $q$-transverse metric on $\tilde{E}$. Then $\eta$ is $q\circ p$-transverse if and only if $q_* \eta$ is $p$-transverse. In that case, $(q\circ p)_* \eta = p_* (q_* \eta)$.
    
\end{lemma}

\subsection{Transversely invariant analytic metrics}

Let $\theta: \G \curvearrowright E$ be an analytic groupoid action with momentum map $q: E \rightarrow M$. Let $\O$ be an orbit of the action. Restricting the action groupoid $\G \ltimes E$ to the orbit $\O$, we get the analytic groupoid
$$\G \ltimes E|_{\O} = \G\ltimes \O = \{(g,e)\in \G \times \O|s(g)= q(e) \}.$$
The groupoid $\G \ltimes \O$ acts linearly on the normal bundle $\N_{\O} = T_\O E/T\O$ as follows. Any $(g,e)\in \G\ltimes \O$ induces a linear isomorphism
\[(g,e): \N_e \rightarrow \N_{ge},\quad (g,e)[v]=\left[\left.\frac{\d}{\d t}\right\vert_{t=0} g(t)e(t)\right],\]
where $g(t)$ and $e(t)$ are curves in $\G$ and $E$ such that $g(0)=g$, $e(0)=e$, $e'(0)=v$ and $s(g(t))= q(e(t))$. One can check that this definition is independent of choice of curves. One obtains a groupoid action, called the \textbf{normal representation}, denoted by $\N_{\O}(\theta)$.

We can dualize the normal representation to define the \textbf{conormal representation}
\[ \N^{*}_{\O}(\theta): \G\ltimes \O \curvearrowright \N^{*}_{\O}. \] Identifying $\N^*_{O} \cong (T\O)^\circ$, this representation is given by:
\[ \N^*_{\O}(\theta)_{(g,e)} (\alpha) = \alpha \circ \N_{\O}(\theta)_{(g^{-1},g e)}. \]

To explicitly express the normal and conormal representations, let $A$ be the Lie algebroid of $\G$ and let $\sigma$ be a right-splitting of the exact sequence
\begin{equation}
\label{sq:spliting_TG}
    \xymatrix{
    0\ar[r] & t^*A \ar[r] & T \G \ar[r]^{s_*} & s^*TM \ar[r] \ar@/^/@{-->}[l]^{\sigma} & 0.}
\end{equation}
Notice that such a splitting $\sigma$ always exists and can be taken to be analytic (use an analytic Riemannian metric on $\G$ as in the proof of Lemma \ref{lem:analytic:p:transv}). The splitting allows us to express the normal representation as:
\begin{equation}\label{eq:normal representation}
    \N_{\O}(\theta)_{(g,e)}([v]) = \left[ \d\theta (\sigma_g(\d_e q(v), v) ) \right].
\end{equation}
On the other hand, the conormal representation can be expressed as:
\begin{equation}\label{eq:conormal representation}
    \langle \N_{\O}^*(\theta)_{(g,e)}(\alpha) , v \rangle = \langle \alpha,  \d\theta (\sigma_{g^{-1}}(d_{ge} q(v), v) \rangle .
\end{equation}
The above expressions do not depend on the choice of splitting $\sigma$.

Now fix an analytic metric $\eta$ on $TE$ and identify the normal bundle $\N_{\O} \cong (T\O) ^ \perp$. The metric $\eta$ induces a dual analytic metric $\eta^{*}$ on the cotangent bundle $T^{*}E$. Restricting the metric $\eta$ to $(T\O)^\perp \cong \N_{O} $ and the metric $\eta^*$ to $(T\O)^\circ \cong \N^*_{O} $, we get analytic metrics on $\N_{\O}$ and $\N_{\O}^*$ respectively.

\begin{definition} Let $\theta: \G \curvearrowright E$ be a groupoid action. We say that a metric $\eta$ on $E$ is \textbf{transversely $\theta$-invariant} if for any orbit $\O$, the normal representation $\N_{\O}(\theta): \G\ltimes \O \curvearrowright \N_{\O}$ acts by isometries or, equivalently, if the conormal representation $\N^{*}_{\O}(\theta): \G\ltimes \O \curvearrowright \N^{*}_{\O}$ acts by isometries.
\end{definition}

Recall that a groupoid $\G$ is called:
    \begin{enumerate}[(i)]
        \item \textbf{proper} if the map $(s,t): \G \rightarrow M \times M$ is proper,
        \item \textbf{s-proper} if the source map $s: \G \rightarrow M$ is proper.
    \end{enumerate}
    
We recall also that a groupoid action is proper if the corresponding action groupoid is proper. The following result, which is proved in \cite{del2018riemannian} in the smooth category, is easily seen to hold also in the analytic setting.

\begin{lemma}
\label{lm:quotient_metric}
Let $\theta:\G \curvearrowright E$ be free and proper groupoid action with quotient map $\pi: E\rightarrow E/\G$. A metric $\eta$ on $E$ is transversely $\theta$-invariant if and only if $\eta$ is $\pi$-transverse. In that case, $E/\G$ has a push-forward metric $\pi_* \eta$, which makes $\pi$ a Riemannian submersion.
\end{lemma}

\subsection{Quasi-actions and analytic metrics}

To construct 2-metrics on a groupoid by averaging one needs the following generalization of a groupoid action.

\begin{definition}
    Let $\G \tto M $ be a Lie groupoid, let $E$ be a manifold with a map $q: E \rightarrow M$. Denote $\G\times_M E = \{ (g,e)\in \G \times E \mid s(g) = q(e) \}$. We say that $\theta: \G \times_{M} E \rightarrow E $ is a \textbf{quasi-action} if $q(\theta(g,e)) = t(g)$, for all $(g,e)\in \G \times_{M} E$.
\end{definition}

We denote a quasi-action by the symbol $\theta : \G \tilde{\curvearrowright} E$.
The quasi-action $\theta$ associates to an arrow $g: y \leftarrow x$ a map $\theta_g: E_x \rightarrow E_y$ between the fibers of the moment map $q$. A groupoid action $\theta: \G \curvearrowright E$ may not admit a lift to a groupoid action on the tangent bundle $TE$, but can always be lifted to a quasi-action with the help of a splitting of \eqref{sq:spliting_TG}.

\begin{definition}
  Let $\sigma: s^* TM \rightarrow T\G$ be a splitting of the  sequence \eqref{sq:spliting_TG} and let $\theta: \G\curvearrowright E$ be a groupoid action with moment map $q: E\rightarrow M$. The \textbf{tangent lift} of $\theta$ is the groupoid quasi-action $T_{\sigma}\theta: \G\ltimes E \tilde{\curvearrowright} TE$, with moment map the projection $p: TE\rightarrow E$, defined by:
    \begin{equation}\label{eq:tangent lift}
      T_{\sigma}\theta_{(g, e)}(v) := \d\theta (\sigma_g(d_e q(v)), v). 
  \end{equation}
    The \textbf{cotangent lift} is the quasi-action $T_{\sigma}^*\theta: \G \ltimes E \tilde{\curvearrowright} T^* E$ defined by:
    \begin{equation}\label{eq:cotangent lift}
      \langle T_{\sigma}^*\theta_{(g, e)}(\alpha) , v \rangle := \langle \alpha,  T_{\sigma}\theta_{(g^{-1}, ge)}(v)\rangle.
  \end{equation}
\end{definition}

\begin{remark}
    In \cite{del2018riemannian} the tangent and cotangent lift were defined using a connection, but for this one only needs a splitting of \eqref{sq:spliting_TG}.
\end{remark}

In the analytic setting, the tangent lift $T_{\sigma}\theta$ and the cotangent lift $T_{\sigma}^*\theta$ are also analytic. The following theorem generalizes a result of \cite{del2018riemannian} to the analytic setting.

\begin{proposition}\label{thm:representation of cotangent lift}
    Let $\theta: \G\curvearrowright E$ be a groupoid action with moment map $q:E\rightarrow M$, let $\O$ be an orbit of the action and fix a splitting of  \eqref{sq:spliting_TG}. Then the conormal bundle $(T\O)^\circ$ is invariant under the cotangent lift quasi-action $T_{\sigma}^*\theta: \G \ltimes E \tilde{\curvearrowright} T^* E$ and the restriction of the quasi-action to  the conormal bundle $(T\O)^\circ$ agrees with the conormal representation: $T_{\sigma}^*\theta \mid_{(T\O)^\circ} = \N_{\O}^*(\theta) $.
\end{proposition}
\begin{proof}
Restricting the action to the orbit 
\[\theta\mid_{\O}: \G\times_M \O \rightarrow \O, \]
and differentiating, we find that
\[\d\theta\mid_{\O}(T\G \times_{TM} T\O )\subset T\O. \]
Hence, for any $v\in T_{e}\O$, we have
\[ T_{\sigma}\theta_{(g, e)}(v) = d\theta (\sigma_g(d_e q(v), v) \in T_{ge}\O. \] 
This shows that $T\O$ is invariant under the quasi-action $T_{\sigma}\theta$. Moreover, by \eqref{eq:cotangent lift}, $(T\O)^\circ$ is also invariant under the cotangent quasi-action  $T_{\sigma}^* \theta$. Upon restriction to $(T\O)^\circ$, the conormal representation \eqref{eq:conormal representation} and cotangent lift \eqref{eq:cotangent lift} coincide.
\end{proof}

\subsection{Averaging}
Henceforth, given a groupoid action $\theta: \G\curvearrowright E$ and a splitting $\sigma$ of \eqref{sq:spliting_TG}, for simplicity, we write
\[ ge:= \theta_{g}(e), \quad g\alpha:=T_{\sigma}^*\theta_{(g,e)}(\alpha).\]

Assume now that the groupoid $\G$ is s-proper and let $\mu$ be an analytic normalized Haar density on $\G$, which always exists -- see Appendix \ref{sec:analytic_haar_density}. Given an analytic metric $\eta$ on $E$, we can average the metric by considering the dual metric, as in \cite{del2018riemannian}, to obtain a transversely invariant metric. Namely, for any $x = q(e)$ and $\alpha, \beta \in T^*_eE$ we set
\begin{equation}
    \label{eq:average}
    (\tilde{\eta})^*_e(\alpha, \beta) := \int_{g\in G(-,x)} \eta^*_{ge}(g\alpha, g\beta) \, \mu^x(g). 
\end{equation}
We call the metric $\tilde{\eta}$ the \textbf{cotangent average} of $\eta$. We have the following analytic version of a result of \cite{del2018riemannian}.

\begin{theorem}
\label{thm:average}
    If $\G$ is an s-proper analytic groupoid, $\theta: \G\curvearrowright E$ is an analytic groupoid action and $\eta$ is an analytic metric, the cotangent average $\tilde{\eta}$ is a transversely $\theta$-invariant, analytic, metric.
\end{theorem}

\begin{proof}
According to Lemma \ref{analyticity}, $\tilde{\eta}$ is analytic, so it suffices to show that for any orbit $\O\subset E$, the conormal representation $\N^*_{\O}(\theta): \G \curvearrowright (T\O)^{\circ}$ acts by isometries with respect to the metric $\tilde{\eta}^*|_{(T\O)^{\circ}\times T\O)^{\circ}}$. According to Theorem \ref{thm:representation of cotangent lift}, for any $\alpha \in (T\O)^\circ$, one has
\[\N_{\O}^*(\theta)_{(g,e)}(\alpha) = g\alpha,\]
so for $\alpha \in (T\O)^\circ$ and any pair of composable arrows we have the usual action identity:
\[ h(g\alpha)=(hg)\alpha.\]
Using this, we find for $\alpha,\beta \in (T\O)^\circ$ that
\begin{align*}
(\tilde{\eta})^*_{ge}(g\alpha, g\beta)
    =&\int_{G(-,y)} \eta^*_{h(ge)}(h(g\alpha), h(g\beta)) \, \mu^y(h)\\
    =&\int_{G(-,y)} \eta^*_{(hg)e}((hg)\alpha, (hg)\beta) \, \mu^y(h)\\
    =& \int_{G(-,y)} \eta^*_{ke}(k\alpha, k\beta) \, \mu^y(k)\\
    =&  (\tilde{\eta})^*_{e}(\alpha,\beta),
\end{align*}
where we set $k=hg$ and used the invariance of the Haar density.
\end{proof}

\subsection{Existence of analytic 2-metrics} We have now everything in place to prove existence of analytic 2-metrics. First, we recall their definition.

Let $\G$ be an analytic Lie groupoid. We denote by $\G^{[k]}$ the space of $k$-tuples of arrows with the same target and by $\G^{(k)}$ the space of $k$-tuples of composable arrows. The groupoid $\G$ acts on $\G^{[k]}$ by left multiplication $g(h_1,\cdots, h_k) = (gh_1, \cdots, gh_k)$. This action is free and proper with quotient map 
\[ \pi^{(k)}: \G^{[k]} \rightarrow \G^{(k-1)} \cong \G^{[k]}/\G,\quad (g_1, \cdots, g_k)\mapsto (g_1^{-1}g_2,\cdots, g_{k-1}^{-1}g_k). \]
For $\G^{(2)}=\G_s\times_t \G$ we have the multiplication map and the two projections on each factor
\[ m,\pi_1,\pi_2:\G^{(2)}\to \G. \]

\begin{definition}
An {\bf analytic 2-metric} in $\G$ is an analytic Riemannian metric $\eta^{(2)}$ on $\G^{(2)}$ which is $m,\pi_1, \pi_2$-transverse and the push-forward metrics along these maps coincide:
\[(\pi_1)_*\eta^{(2)}=(\pi_2)_*\eta^{(2)}  = m_*\eta^{(2)}.\]
\end{definition}

\begin{theorem}
\label{thm:analytic:2:metric}
Every analytic s-proper groupoid admits an analytic 2-metric.
\end{theorem}

\begin{proof}
The proof is entirely similar to the smooth case -- see \cite{del2018riemannian}.
By Lemma \ref{lem:analytic:p:transv}, we can choose an analytic $t$-transverse metric $\eta$ on $\G$. It induces a k-fold metric $\eta^{[k]}$ on $\G^{[k]}$ which is transverse to the the map $\G^{[k]}\to M$. By Theorem \ref{thm:average} its cotangent average is a metric $\Tilde{\eta}^{[k]}$ on $\G^{[k]}$ which is transversely invariant for the left $\G$-action on $\G^{[k]}$. Hence, by Lemma \ref{lm:quotient_metric}, it induces a push-forward metric $\tilde{\eta}^{(k-1)} = \pi^{(k)}_{*}\tilde{\eta}^{[k]}$ on $\G^{(k-1)}$. We claim that $\tilde{\eta}^{(2)}$ is an analytic 2-metric on $\G^{(2)}$.

For that we just need to observe that we have a commutative diagram
\[
\xymatrix@1{
\cdots \,
\G^{[3]} \ar@<0.5pc>[r]^{p_3} \ar[r]|{p_2} \ar@<-0.5pc>[r]_{p_1} \ar[d]_{\pi^{(3)}}&
\G^{[2]} \ar@<0.25pc>[r]^{\pr_1} \ar@<-0.25pc>[r]_{\pr_2} \ar[d]_{\pi^{(2)}}& 
\G\ar[d]^{s=\pi^{(1)}}\\
\cdots \,
\G^{(2)} \ar@<0.5pc>[r]^{\pi_1} \ar[r]|{m} \ar@<-0.5pc>[r]_{\pi_2} &
\G \ar@<0.25pc>[r]^{t} \ar@<-0.25pc>[r]_{s} & 
M}
\]
where all the arrows are Riemannian submersions. Here, $p_i:\G^{[3]}\to \G^{[2]}$ is the map that omits the arrow $g_i$:
\[ p_1(g_1,g_2,g_3)=(g_2,g_3),\quad p_2(g_1,g_2,g_3)=(g_1,g_3),\quad p_3(g_1,g_2,g_3)=(g_1,g_2).\]
Since the diagram is commutative, we have 
\[ (\pi_1)_*\tilde{\eta}^{(2)}=(\pi_2)_*\tilde{\eta}^{(2)}  = m_*\tilde{\eta}^{(2)} = \pi^{(2)}_{*}\tilde{\eta}^{[2]}. \]
\end{proof}

\section{Analytic linearization of s-proper groupoids}
Let $\G \tto M$ be a Lie groupoid. We say that a submanifold $S\subset M$ is an \textbf{invariant submanifold} if $S$ is a union of $\G$-orbits. Given such invariant submanifold, we let
\[ \G_{S} := s^{-1}(S),\quad  \N_{\G_S} := T_{\G_S}\G/T\G_S,\quad \N_S:=T_S M/TS.\] 
The normal bundle $\N_{\G_S}$ has a groupoid structure over $\N_S$ with structure maps induced by the differential of structure maps of $\G$, i.e., by $\d s,\d t,\d m,\d\iota$. We call the groupoid $\N_{\G_S}\tto \N_S$ the \textbf{linear model} of $\G$ around $S$. 

There are various other ways of describing the groupoid structure of the linear model. For example, the normal bundle $\N_{\G_{S}}$ coincides with the pull-back of $\N_S$ along the source map
\[ \N_{\G_{S}} = \left\{ (g,v)\in \G_S\times \N_S| s(g)=\pi(v) \right\}.\]
The groupoid $\G_S$ acts on $\N_S$ by the normal representation and the groupoid $\N_{\G_S}\tto \N_S$ is just the action groupoid $\G_S \ltimes \N_S\tto \N_S$.

\begin{definition}
Let $\G\tto M$ be an analytic groupoid and let $S\subset M$ be an invariant submanifold. We say that:
\begin{enumerate}[(i)]
    \item $\G$ is {\bf analytically linearizable} around $S$ if there exist open neighborhoods of $S$, $U\subset M$ and $V\subset \N_{S}$, and an analytic groupoid isomorphism
   \[ \G|_{U}\cong \N_{\G_{S}}|_{V}\]
   which is identity on $\G_{S}$.
   \item $\G$ is {\bf analytically invariantly linearizable} around $S$ if the open neighborhoods $U$ and $V$ can be chosen to be invariant.
\end{enumerate}
\end{definition}

The fundamental observation made in \cite{del2018riemannian} is that the exponential map of a 2-metric yields a linearization map. This remains true in the analytic category.

\begin{theorem}[\cite{del2018riemannian}]
    Let $\G\tto M$ be an analytic groupoid and let $S\subset M$ be an invariant submanifold. If $\G$ admits an analytic 2-metric then it is analytically linearizable around $S$.
\end{theorem}
    \begin{proof}
        An analytic 2-metric $\eta^{(2)}$ on $\G^{(2)}$ induces an analytic metric $\eta^{(1)}=m_{*}\eta^{(2)}$ on $\G$ and their exponential maps gives an analytic linearization of $\G$ around $S$. Namely, there exist open neighborhoods $U\subset M$ and $V\subset \N_{S}$ of $S$ where the exponential maps restrict to a groupoid isomorphism:
    \[
    \xymatrix{
    \N_{\G_{S}}|_{V} \ar[r]^{\exp^{(2)}} \ar@<0.2pc>[d] \ar@<-0.2pc>[d] & \G|_{U} \ar@<0.2pc>[d] \ar@<-0.2pc>[d] \\
    V \ar[r]^{\exp^{(1)}} & U
    }
    \]
    \end{proof}
    
If we further assume that $\G$ is s-proper, then $U$ and $V$ in the previous proof can be chosen to be invariant so that $\G$ is analytically invariantly linearizable around $S$. Therefore, as a consequence of Theorem \ref{thm:analytic:2:metric} we deduce our first main result.

\begin{theorem}\label{thm:analytic_linearization}
    Let $\G\tto M$ be an analytic s-proper Lie groupoid and let $S\subset M$ be an invariant submanifold. Then $\G$ can be analytically invariantly linearizable around $S$.
\end{theorem}

This result can be applied to obtain analytic versions of known linearization results in the smooth category. For example, in Poisson geometry, one obtains the following analytic version of \cite[Theorem 8.6]{FM24}:

\begin{corollary}
Let $(M,\pi)$ be an analytic Poisson manifold and $S\subset M$ a Poisson submanifold. If $T^*_SM$ is integrable by a compact, Hausdorff, Lie groupoid whose source fibers have trivial 2nd de Rham cohomology, then $(M,\pi)$ is analytically linearizable around $S$.
\end{corollary}

\begin{proof}
Theorem 8.2 in \cite{FM24} states that a Poisson manifold $(M,\pi)$ is linearizable around a Poisson submanifold $S$ if and only if the algebroid $T^*M$ is linearizable around $S$. The proof in \cite{FM24} holds in the analytic category, so all we need to show is that the algebroid $T^*M$ is analytically linearizable around $S$. If we show that $(M,\pi)$ is integrable by an analytic s-proper Lie groupoid on an invariant neighborhood $U$ of $S$, then Theorem \ref{thm:analytic_linearization} shows that it is analytically linearizable around $S$, so its algebroid will be analytically linearizable around $S$, and the proof will be completed.

Corollary 8.9 in \cite{FM24} shows that the assumptions made in the statement imply that there is an invariant neighborhood $U$ of $S$ such that $(U,\pi|_U)$ is integrable by an s-proper Lie groupoid. But this groupoid has a unique compatible real analytic structure, since it integrates an analytic algebroid: this follows either by observing that the exponential map gives a real analytic structure in a neighborhood of the identity section and any such neighborhood generates the whole groupoid (Proposition 13.7 in \cite{CFM21}) or, alternatively, by using the result of \cite{Martinez20} which shows that every proper groupod has a compatible analytic structure.
\end{proof}

\section{Holomorphic extensions of Lie algebroids}
\label{sec:algebroids:holomorphic:extensions}

Given a real Lie algebra $\mathfrak{g}$, one can complexify the underlying vector space and extend the Lie algebra structure linearly to obtain the complex Lie algebra $\mathfrak{g}_{\mathbb{C}}$. One calls $\mathfrak{g}$ a real form of $\mathfrak{g}_{\mathbb{C}}$ and $\mathfrak{g}_{\mathbb{C}}$ the complexification of $\mathfrak{g}$. A generalization of these notions to analytic Lie algebroids was proposed by Weinstein in \cite{weinstein2007integration}. We will start by recalling this notion and then we discuss some of its properties. Before, we collect some useful facts about totally real submanifolds.

\subsection{Totally real submanifolds}

Let $(X,J)$ be a complex manifold and let $M\subset X$ be an embedded real analytic submanifold of $X$. We recall that $M$ is called a \textbf{totally real submanifold} of $X$ if  
\[TM\oplus J(TM)=TX|_M. \]

A classical result of Bruhat and Whitney \cite{whitney1959quelques} states that every real analytic manifold $M$ admits a totally real closed embedding 
\[ \iota:M\to M_\C \] 
into some complex manifold $M_\C$. The manifold $M_\C$ is not unique, but the germ of $M$ in $M_\C$ is: given two totally real embeddings $\iota_k:M\to X_k$, there exists open neighborhoods $V_k\subset X_k$ of $M$ and a holomorphic isomorphism $\phi:V_1\to V_2$ such that $\phi\circ \iota_1=\iota_2$.

\begin{remark}
\label{rem:extending:cover}
    A useful fact is a classical result of Grauert \cite{Grauert58} showing that one can take $M_\C$ to be an open neighborhood of the zero section $\iota:M\to TM$. For example, it follows that if one is given an open cover $\{U_i\}$ of $M$ with $U_i\cap U_j$ connected, then one can assume that $M_\C$, possibly after shrinking, also admits a cover $\{U_i^*\}$ such that 
    \[ U^*_i\cap M= U_i, \quad \overline{U^*_i}\cap M= \overline{U_i},\]
    and the intersections $U^*_i\cap U^*_j$ are connected. Here, and henceforth, an overline over a set denotes its closure, unless otherwise stated.
\end{remark}

The following result is well-known and not hard to prove. We state it for future reference. 

\begin{lemma}\label{lem:totally real}
Let $M_\C$ be a complex manifold and let $M\subset M_\C$ be an embedded real analytic submanifold of $M_\C$. The following are equivalent:
\begin{enumerate}[(i)]
\item $M$ is a totally real submanifold of $M_\C$
    \item\label{local normal form} For any $x\in M$, there is a  chart $U$ for $M_\C$ centered at $x$, such that 
$$M\cap U=\{(x_{1},\dots,x_{n})\in \mathbb{C}^n \mid x_{i}\in \mathbb{R}\quad i=1,\dots,n\}$$
\item There is an open neighborhood $V\subset M_\C$ containing $M$ and an anti-holomorphic involution $\sigma: V\xrightarrow{} V$ with $\mathrm{Fix}(\sigma)=M$
\end{enumerate}
Moreover, if $Y$ is a holomorphic manifold and $\phi:M\to Y$ is a real analytic map, there exists an open $V\subset M_\C$ and a holomorphic map $\phi^*:V\to Y$ such that $\phi=\phi^*\circ\iota$.
\end{lemma}

Notice that the germ of the holomorphic extension $\phi^*$, in the last part of the statement, is unique: if $\phi^*_k:V_k\to Y$, $k=1,2$, are two holomorphic extensions of $\phi:M\to Y$, then there exists an open $V\subset V_1\cap V_2$ containing $M$ such that $\phi^*_1|_V=\phi^*_2|_V$.



\subsection{Definition of holomorphic extension} 
Let $(A, \rho,[\cdot,\cdot])$ be an analytic Lie algebroid, so $A$ is an analytic vector bundle over the analytic manifold $M$, $\rho$ is an analytic map, and the bracket $[\cdot,\cdot ]$ takes analytic sections to analytic sections. It follows from \cite{crainic2003integrability} that any Lie groupoid integrating an analytic algebroid $A$ is automatically analytic.

\begin{definition}
\label{def:holomorphic:extension:algebroid}
Let $(A, \rho,[\cdot,\cdot])$ be an analytic Lie algebroid over an analytic manifold $M$. A holomorphic algebroid $(A_{\C},\rho_{\C},[\cdot,\cdot]_{\C})$ over a complex manifold $M_{\C}$ is called a \textbf{holomorphic extension} of $A$ if the following hold:
    \begin{enumerate}[(i)]
        \item $M$ is a totally real submanifold of $M_{\C}$;
        \item $A_{\C}|_{M}= A\oplus J(A)$, where $J$ is the complex structure on $A_{\C}$;
        \item $\rho_\C$ is an extension of $\rho$, i.e., the following diagram is commutative
        \begin{equation}\label{diag: algebroid_complexify}
        \vcenter{\xymatrix{
        A \ar[r] \ar[d]_{\rho}
        & A_{\C} \ar[d]^{\rho_{\C}} \\
        TM \ar[r] &  T M_{\C}}}
       \end{equation}
\item $[\cdot,\cdot]_\C$ an extension of $[\cdot,\cdot]$, i.e., for any $x\in M$, any local analytic sections $\alpha$ and $\beta$ of $A$ with local holomorphic extensions $\Tilde{\alpha}$ and $\Tilde{\beta}$
\[[\alpha,\beta](x)=[\Tilde{\alpha},\Tilde{\beta}]_\C(x).  \]
    \end{enumerate}
\end{definition}

\begin{remark}
\label{rem:real:form}
    Under the assumptions of this definition, one can also call $A\Rightarrow M$ a \textbf{real form} of $A_\C\Rightarrow M_\C$. However, at the groupoid level these two notions become distinct. For this reason, we will not use this term.
\end{remark} 

Lie algebroids often arise from some geometric structure. In this case, their holomorphic extensions arises from a holomorphic extension of the geometric structure. Here is a simple example.

\begin{example}
\label{ex:Poisson:holomorphic}
    Let $(\mathbb{R}^n, \pi)$ be an analytic Poisson manifold. Write
   \[\pi=\sum_{i,j=1}^{n}\pi_{ij}(x)\frac{\partial}{\partial{x_{i}}}\wedge\frac{\partial}{\partial{x_{j}}}\]
   This gives rise to an analytic Lie algebroid $(T^*\R^n, \pi^{\#},[\cdot,\cdot]_{\pi})$.
   
   Now consider the standard embedding $\mathbb{R}^n\subset \mathbb{C}^n$. Since the components $\pi_{ij}$ are analytic, we can extend them to get holomorphic functions $\tilde{\pi}_{ij}$ in a connected open neighborhood $U \subset \C^n$ of $\mathbb{R}^n$. This gives rise to a holomorphic bivector
      \[\pi_{\C}=\sum_{i,j=1}^{n}\tilde{\pi}_{ij}(z)\frac{\partial}{\partial{z_{i}}}\wedge\frac{\partial}{\partial{z_{j}}}.\]
    The Schouten bracket of $\pi_\C$ is a holomorphic extension of the Schouten bracket of $\pi$, so it vanishes. Hence, $\pi_\C$ is a holomorphic Poisson structure and one checks easily that its cotangent algebroid $(T^*U, \pi_{\C}^{\#},[\cdot,\cdot]_{\pi_{\C}})$ is a holomorphic extension of $(T^*\R^n, \pi^{\#},[\cdot,\cdot]_{\pi})$.

    More generally, any analytic Poisson manifold $(M,\pi)$ has a holomorphic extension $(M_\C,\pi_\C)$ whose holomorphic cotangent algebroid $(T^*M_\C, \pi_{\C}^{\#},[\cdot,\cdot]_{\pi_{\C}})$ is a holomorphic extension of $(T^*M, \pi^{\#},[\cdot,\cdot]_{\pi})$.
\end{example}

\subsection{Existence} The holomorphic extension of an algebroid $A$ is not unique, but its germ is unique: if $A_1$ and $A_2$ are two holomorphic extensions of $A$, then they are isomorphic when restricted to open neighborhoods of $M$. This follows from uniqueness of holomorphic extensions. Moreover, it is not hard to show that holomorphic extensions always exist.

\begin{proposition}
     Every analytic algebroid $A\Rightarrow M$ admits a holomorphic extension.
\end{proposition}
   
\begin{proof}
    According to \cite{whitney1959quelques,shutrick1958complex,lewis2023geometric}, there exists a holomorphic vector bundle $A_{\C}\rightarrow M_\C$ such that $M$ is a totally real submanifold of $M_\C$ and $A_{\C}|_{M}= A\oplus J(A)$. Locally we can always choose trivializations and charts such that the embedding $i:A\rightarrow A_\C$ looks like the standard embedding of $\R^n\times\R^r$ in $\C^n\times\C^r$, where $r=\rank_\R A=\rank_\C A_\C$: 
    \[
    \xymatrix{
    A|_{U\cap M} \ar[r]^i \ar[d]_{\simeq}
    & A_{\C}|_{U} \ar[d]^{\simeq} \\
    \mathbb{R}^n\times  \mathbb{R}^r\ar[r]
    &\mathbb{C}^n\times  \mathbb{C}^r
    }\]
Then locally $\rho$ can be expressed as a function $\R^n \rightarrow M_{n\times r}(\mathbb{R})$, which locally can be extended to a holomorphic function $\C^n\rightarrow M_{n\times r}(\mathbb{C})$. Possibly by shrinking $M_\C$ to an open neighborhood of $M$, we can glue those extensions to get a holomorphic anchor map $\rho_\C:A_{\C}\rightarrow T M_{\C} $. Similarly, we can define the bracket $[\cdot,\cdot]_\C$ by holomorphically extending $[\cdot,\cdot]$.
\end{proof}

\subsection{Anti-holomorphic involutions}

Let $A_\C\Rightarrow M_\C$ be a holomorphic Lie algebroid. Forgetting the complex structure on the vector bundle $A_\C$, we can view $A_\C$ as a real vector bundle and the anchor map $\rho_\C: A_\C \rightarrow TM_\C$ as a real bundle map. According to \cite{laurent2008holomorphic}, there exists a unique real Lie algebroid structure on $A_\C$ with the same anchor map and such that the Lie bracket when restricted to holomorphic sections coincides with the Lie bracket of $A_\C$. We denote the real Lie algebroid structure by $(A_R, \rho_R, [\cdot, \cdot]_R)$.

\begin{proposition}\label{prop: involution_of_lie_algebroid}
Let $(A_{\C},\rho_{\C},[\cdot,\cdot]_{\C})\Rightarrow M_{\C}$ be a holomorphic Lie algebroid. Let $A_{R}$ be the underlying real Lie algebroid. Let $A\rightarrow M$ be an analytic subbundle of $A_{\C}\rightarrow M_{\C}$. Then $A_{\C}$ is a holomorphic extension of $A$ if and only if there exists an open neighborhood $V\subset M_\C$ containing $M$ and an anti-holomorphic Lie algebroid involution $\sigma: A_{R}|_{V}\rightarrow A_{R}|_{V}$ with $\mathrm{Fix}(\sigma)=A$.
\end{proposition}

\begin{proof}
    Given a holomorphic Lie algebroid $(A_{\C},\rho_{\C},[\cdot,\cdot]_{\C})\Rightarrow M_{\C}$ with an anti-holomorphic Lie algebroid involution $\sigma: A_{R}|_{V}\rightarrow A_{R}|_{V}$ covering an anti-holomorphic involution $\lambda:V\to V$ on some open set $V\subset M_\C$, we set 
    \[ M:=\mathrm{Fix}(\lambda)\text{ and } A:=\mathrm{Fix}(\sigma).\] 
    Then one checks easily that $A\Rightarrow M$ is a (real) Lie algebroid, $M_\C$ is a complexification of $M$, $A_\C|_M=A\oplus J(A)$, and the bracket and anchors of $A|_\C$ are extensions of the bracket and anchors of $A$.
    
    Conversely, assume $A_{\C}$ is a holomorphic extension of $A$. Locally we can always choose trivializations and charts such that the embedding $i:A\rightarrow A_\C$ looks like the standard embedding of $\R^n$ in $\C^n$: 
    \[
    \xymatrix{
    A|_{U\cap M} \ar[r]^i \ar[d]_{\simeq}
    & A_{\C}|_{U} \ar[d]^{\simeq} \\
    \mathbb{R}^n\times  \mathbb{R}^r\ar[r]
    &\mathbb{C}^n\times  \mathbb{C}^r
    }\]
    Write $(z_k,u_l)$ for coordinates in $\C^n\times\C^r$. Let $\sigma_{U}$ be defined by requiring 
    \begin{equation}
        \label{eq:involution}
        \sigma_{U}(z_{1},\cdots,z_{n},u_1,\cdots,u_r)=(\Bar{z}_1,\cdots,\Bar{z}_n,\Bar{u}_1,\cdots,\Bar{u}_r).
    \end{equation}
    By analytic continuation, we can glue such $\sigma_{U}$ to get an anti-holomorphic involution $\sigma$ on a neighborhood of $M$ with $\mathrm{Fix}(\sigma)=A$. We claim that $\sigma$ is a Lie algebroid morphism.

 To prove our claim, let $\lambda: V \rightarrow V$ be the anti-holomorphic involution covered by $\sigma$. Then $\d\lambda\circ \rho_R \circ \sigma: A_R|_V \rightarrow TV $ is a holomorphic map which when restricted to $A$ coincides with $\rho_\C$. By analytic continuation, 
\[ \d\lambda\circ \rho_R \circ \sigma = \rho_\C,\]
or equivalently
\[ \d\lambda\circ \rho_R  = \rho_R \circ \sigma,\]
so $\sigma$ preserves anchors. To check that it preserves brackets, observe that any smooth section of $A$ can be locally written as linear combination $\sum_{k}f_k e_k$, where $e_k$ is a holomorphic section and $f_k$ is a complex-valued smooth function. To finish proving our claim, it suffices to show that for any smooth complex valued functions $f$ and $g$ we have
\[[\sigma_*(fe_i),\sigma_*(ge_j)]_R =\sigma_*[fe_i, g e_j]_R, \]
where $\sigma_*(e) = \sigma\circ e \circ \lambda^{-1}$. This follows by the following straightforward computation
\begin{align*}
    \sigma_*[fe_i, g e_j]_R(z) &= \overline{[fe_i, g e_j]_R(\bar{z})} \\
     =& \overline{f(\bar{z})g(\bar{z})}[e_i, e_j] +\overline{\rho_R(e_i)(g)(\bar{z})}\cdot \overline{f(\bar{z})}e_j - \overline{\rho_R(e_j)(f)(\bar{z})}\cdot \overline{g(\bar{z})}e_i \\
    =& \overline{f(\bar{z})g(\bar{z})}[e_i, e_j] +\rho_R(e_i)(\overline{g(\bar{z})})\cdot \overline{f(\bar{z})}e_j - \rho_R(e_j)(\overline{f(\bar{z})})\cdot \overline{g(\bar{z})}e_i \\
    =&[\overline{f(\bar{z})}e_i,\overline{g(\bar{z})}e_j ]_R\\
    =& [\sigma_*(fe_i),\sigma_*(ge_j)]_R,
\end{align*}
where we used the local definition of $\sigma$ given by \eqref{eq:involution}.
\end{proof}

\subsection{Universal property}

The holomorphic extension of an algebroid satisfies the following ``universal property''.

\begin{proposition}
    Let $A_\C\Rightarrow M_\C$ be a holomorphic extension of a real analytic Lie algebroid $A\Rightarrow M$ with canonical embedding $\iota:A\to A_\C$. Then
    \begin{enumerate}[(i)]
    \item for any holomorphic algebroid $B\Rightarrow X$ and real analytic morphism $\phi:A\to B$ there exists a neighborhood $V$ of $M$ and a holomorphic morphism $\phi^*:A_\C|_V\to B$ with $\phi=\phi^*\circ\iota$;
    \item if $\psi:A_\C|_V\to B$ is a morphism of holomorphic Lie algebroids, defined on some neighborhood $V$ of $M$, and $\phi=\psi\circ\iota$, then $\psi=\phi^*$ on some neighborhood $V'\subset V$ of $M$ in $M_\C$.
\end{enumerate}
\end{proposition}

\begin{proof}
    In order to prove (i), one can first find local holomorphic extensions of $\phi$ by working in local coordinates, as in the proof of Proposition \ref{prop: involution_of_lie_algebroid}. By uniqueness of analytic continuation, there is an open neighborhood $V$ of $M$ where these local extensions glue to a well-defined  holomorphic morphism $\phi^*:A_\C|_V\to B$ with $\phi=\phi^*\circ\iota$. Property (ii) follows from the uniqueness of holomorphic extensions of maps.    
\end{proof}

\section{Holomorphic extensions of Lie groupoids}
\label{sec:groupoids:holomorphic:extensions}

The notion of holomorphic extension of a groupoid is more subtle than the one of an algebroid. Some inspiration can be drawn from the case of Lie groups.

Recall that given a real Lie group $G$ its {\bf universal complexification} is a complex Lie group $G_\C$ together with a continuous morphism $\iota:G\to G_\C$ satisfying the following universal property: for every complex Lie group $H$ and every continuous morphism $\phi:G\to H$ there exists a unique holomorphic morphism $\phi^*:G_\C\to H$ such that $\phi=\phi^*\circ\iota$. Every Lie group admits a complexification which is unique, up to isomorphism. The image of $\iota:G\to G_\C$ is a closed, totally real, subgroup of $G_\C$. See, e.g., \cite{heinzner1994equivariantExt} for the proofs of these facts.

Notice that the morphism $\iota:G\to G_\C$ may have non-trivial, even non-discrete, kernel. For example, one can let $G=(\widetilde{\mathrm{SL}_2}\times \mathbb{S}^1)/C$, where $C$ is a discrete normal subgroup whose projection in $\mathbb{S}^1$ is dense. Hence, the following peculiarities with the notion of complexification can occur:
\begin{itemize}
    \item a Lie group $G$ may fail to be a totally real submanifold of its complexification;
    \item the Lie algebra of the complexification $G_\C$ may not be isomorphic to the complexification of the Lie algebra of $G$.
\end{itemize}
These issues disappear when $\iota:G\to G_\C$ has trivial kernel, in which case one calls $G_\C$ a {\bf holomorphic extension} of $G$. The holomorphic extension of $G$, if it exists, is unique (up to isomorphism). Several classes of Lie groups admit holomorphic extensions, such as compact Lie groups, linear groups and solvable groups (see, e.g., \cite{heinzner1994equivariantExt,OV1994}).


\subsection{Definition of holomorphic extension}
In the sequel, given a Lie groupoid $\G\tto M$ and a subgroupoid $\H\tto N$, by a \emph{groupoid neighborhood} of $\H$ in $\G$ we mean a pair of open subsets $N\subset V\subset M$, $\H\subset\V\subset \G$ such that $\V\tto V$ is a subgroupoid of $\G\tto M$. When $\V=\G|_V$ we call $\V\tto V$ a \emph{full groupoid neighborhood}. When $\G\tto M$ is a proper Lie groupoid every groupoid neighborhood always contains a full groupoid neighborhood (see, Lemma 5.1.3 in \cite{del2018riemannian}).

\begin{definition}
\label{def:holomorphic:extension:groupoid}
A holomorphic groupoid $\G_{\C} \tto M_{\C}$ is called a \textbf{holomorphic extension} of a real analytic Lie groupoid $\G \tto M$ if there exists a real analytic morphism $\iota:\G\to\G_\C$ such that:
\begin{enumerate}[(C1)]
    \item $\iota:\G\to\G_\C$ is a totally real, closed, embedding;
    \item for any holomorphic groupoid $\H$ and real analytic morphism $\phi:\G\to\H$ there exists a groupoid neighborhood $\V$ of $\iota(G)$ in $\G_\C$ and a holomorphic morphism $\phi^*:\V\to \H$ with $\phi=\phi^*\circ\iota$;
    \item if $\psi:\V\to \H$ is a morphism of holomorphic groupoids on some groupoid neighborhood $\V$ of $\iota(\G)$ in $\G_\C$, and $\phi=\psi\circ\iota$, then $\psi=\phi^*$ on some groupoid neighborhood $\V'\subset \V$ of $\iota(\G)$ in $\G_\C$.
\end{enumerate}
A \textbf{full holomorphic extension} is a holomorphic extension where in (C2) the groupoid neighborhood $\V$ can be taken to be full, and where in (C3) whenever $\V$ is full, one can take $\V'$ also to be full.
\end{definition}

An analytic Lie groupoid may not admit a holomorphic extension (as we discussed above, this is already true for Lie groups).  However, it follows from properties (C2) and (C3) in the definition that the germ of a holomorphic extension of $\G$, if it exists, is unique.

Given a Lie groupoid $\G$ we denote its Lie algebroid by $\Lie(\G)$. If $\G$ is a real analytic (respectively, holomorphic) Lie groupoid then $\Lie(\G)$ is a real analytic (respectively, holomorphic) Lie algebroid.

\begin{proposition}\label{prop: complexify}
    Let $\mathcal{G}_{\mathbb{C}}$ be a holomorphic extension of $\G$. Then $\Lie(\G_{\C})$ is a holomorphic extension of $\Lie(\G)$.     
\end{proposition}

\begin{proof}
Since $\iota:\G\to \G_\C$ is totally real, one has
\[ T\G_\C=T\G\oplus J(T\G). \]
Using that
\[ T_M\G=TM\oplus\ker(\d s),\quad
T_{M_\C}\G_\C=TM_\C\oplus\ker(\d s_\C), \]
and the fact that $M_\C$ and the source fibers are complex submanifolds of $\G_\C$, 
one concludes that
\[ TM_\C=TM\oplus J(TM), \quad \ker(\d_M s_\C)=\ker(\d_M s)\oplus J(\ker\d_M s). \]
This shows that the vector bundle $A_\C=\Lie(\G_\C):=\ker(\d_M s_\C)$ is a  holomorphic extension of the vector bundle $A=\Lie(\G):=\ker(\d_M s)$.

Since $\iota:\G\to \G_\C$ is a groupoid morphism, it follows that $i_*\circ\rho=\rho\circ i_*$, so $\rho_\C$ is an extension of $\rho$. On the other hand, let $\alpha,\beta \in  \Gamma(A)$ be analytic sections of $A$ and denote by $\overrightarrow{\alpha}$ and $\overrightarrow{\beta}$ the corresponding right-invariant vector fields in $\G$ so that
\[ [\overrightarrow{\alpha},\overrightarrow{\beta}]=\overrightarrow{[\alpha,\beta]}_A. \]
Extending $\alpha$ and $\beta$ to holomorphic sections $\tilde{\alpha}$ and $\tilde{\beta}$ of $A_\C$, the Lie bracket of the corresponding right-invariant vector fields in $\G_\C$
\[
[\overrightarrow{\tilde{\alpha}}, \overrightarrow{\tilde{\beta}}]=
\overrightarrow{[\tilde{\alpha},\tilde{\beta}]}_{A_\C}
\]
is a holomorphic extension of $[\overrightarrow{\alpha}, \overrightarrow{\beta}]$. It follows that $[\tilde{\alpha}, \tilde{\beta}]_\C$ is a holomorphic extension of $[\alpha,\beta]$.

\end{proof}

\begin{remark}
    The previous proof only used property (C1) in the definition of holomorphic extension.
\end{remark}

\subsection{Examples}
\label{sec:examples:groupois}
We discuss several classes of groupoids which admit holomorphic extensions.

\subsubsection{Unit groupoid}
    Let $M$ be an analytic manifold and consider the unit groupoid $M\tto M$. If $M_\C$ is any complexification of $M$, the unit groupoid $M_\C\tto M_\C$ is a full holomorphic extension of $M\tto M$ (for the canonical embedding $M\hookrightarrow M_\C$). Indeed the properties in the definition follow from Lemma \ref{lem:totally real} and the fact that any groupoid morphism from the unit groupoid $M\tto M$ to some other groupoid $\H\tto N$ amounts simply to a map $M\to N$.
    In other words, the complexifications $M$ are in 1:1 correspondence with the holomorphic extensions of the unit groupoid $M\tto M$.

\subsubsection{Pair groupoid}
    Let $M$ be a real analytic manifold and let $M_\C$ be a complexification of $M$, with totally real embedding $i:M\to M_\C$. The pair groupoid $M_\C\times M_\C\tto M_\C$ is a holomorphic extension of the real analytic groupoid $M\times M\tto M$ for the groupoid embedding 
    \[ \iota:=(i,i):M\times M\to M_\C\times M_\C. \]
    To see this, let $\H\tto X$ be a holomorphic groupoid and $\phi:M\times M\to \H$ a real analytic morphism. Any such morphism takes the form
    \[ \phi(p,q)=\psi(p)\psi(q)^{-1}, \]
    where $\psi:M\to s_\H^{-1}(x_0)$ is an analytic map. This map has a unique holomorphic extension $\psi^*:V\to s_\H^{-1}(x_0)$ defined in some open $V\subset M_\C$ containing $M$. Hence, we can define a holomorphic morphism $\phi^*:V\times V\to \H$ extending $\phi$ by setting
    \[ \phi^*(p,q):=\psi^*(p)\psi^*(q)^{-1}. \]
    Since $\V:=V\times V\tto V$ is a full groupoid neighborhood of 
    $M\times M$ in $M_\C\times M_\C$, this shows that (C2) holds. Also, one checks easily that (C3) holds due to the uniqueness of holomorphic extensions of maps.


\subsubsection{Transitive proper groupoids}    
\label{sec:complexify:proper groupoid}
    Generalizing both compact Lie groups and pair groupoids, let $\G\tto M$ be an analytic transitive proper groupoid. Then $\G$ is isomorphic to the gauge groupoid of an analytic principal bundle $P\to M$ with structure group a compact Lie group $G$:
\[ \G\simeq P\times_G P. \]
For such principal bundles we have the following result (for a proof see also Appendix \ref{sec:covers}).

\begin{theorem}[\cite{heinzner1994equivariantExt}]
\label{thm:complexify_principal_bundle}
    Let $G$ be a compact group and $P\to M$ be an analytic principal $G$-bundle. Let $G_\C$ be the holomorphic extension of $G$. There exists a holomorphic principal $G_\C$-bundle $P_\C\to M_\C$ such that $P$ is a totally real submanifold of $P_\C$ and the embedding $P\rightarrow P_\C$ is $G$-equivariant. 
\end{theorem}

Using this result, we claim that the transitive groupoid $\G=P\times_G P$ has full holomorphic extension the holomorphic transitive groupoid 
\[ \G_\C=P_\C\times_{G_\C} P_\C, \]
with the obvious embedding $\iota:\G\to\G_\C$. Property (C1) clearly holds. To check the other two properties, we use the following.

\begin{lemma}
Given a holomorphic groupoid $\H\tto X$, every analytic morphism $\phi:\G\to \H$ can be written in the form
\[ \phi([p_1,p_2])=\psi(p_1)\psi(p_2)^{-1}, \]
for an analytic map $\psi:P\to s^{-1}(y_0)$ satisfying:
\begin{enumerate}[(a)]
\item $\psi$ is $G$-equivariant relative to a Lie group morphism $\rho:G\to \H_{y_0}$;
\item $\psi$ covers a map $\psi_0:M\to X$.
\end{enumerate} 
\end{lemma}

\begin{proof}
    Fix any $x_0=p_0 G\in M=P/G$ so we can identify the source fiber and isotropy group of $\G=P\times_G P$ at $x_0$ with, respectively, $P$ and $G$, via the embeddings
    \begin{align*}
        i_P:P&\hookrightarrow P\times_G P,\quad p\mapsto [p,p_0],\\ 
        i_G:G&\hookrightarrow P\times_G P,\quad g\mapsto [p_0g,p_0].
    \end{align*}
    If $1_{y_0}=\phi(1_{x_0})$, we define $\psi:P\to s^{-1}(y_0)$ and $\rho:G\to \H_{y_0}$ as the composition of $\psi$ with these embeddings:
    \[ \psi:=\phi\circ i_P, \quad \rho:=\phi\circ i_G. \]
    Then, using that $\phi$ is a groupoid morphism, we find
    \[ \phi([p_1,p_2])=\phi([p_1,p_0]\cdot[p_0,p_2])=\psi(p_1)\psi(p_2)^{-1}. \]
    One also checks easily that (a) and (b) hold.
\end{proof}

In this case, it follows from \cite{heinzner1994equivariantExt} that $\psi$ has a holomorphic extension $\psi^*:\widehat{P}\to s^{-1}(x)$ defined in some $G_\C$-equivariant neighborhood $\widehat{P}$ of $P$ in $P_\C$. This extension is $G_\C$-equivariant relative to the complexification of the morphism $\rho:G\to \H_x$. 

We can then define a holomorphic groupoid morphism
\[ \phi^*:\widehat{P}\times_{G_\C}\widehat{P}\to \H, \quad  \phi^*([p_1,p_2]):=\psi^*(p_1)\psi^*(p_2)^{-1}. \]
Notice that $\widehat{P}\times_{G_\C}\widehat{P}\tto \widehat{P}/G_\C$ is a full groupoid neighborhood of $\G$ in $\G_\C$, and properties (C2) and (C3) now follow.

\subsubsection{Bundle of compact Lie groups} Let $\G\to M$ be a bundle of compact, connected, Lie groups. Such bundles are always locally trivial and it follows that there is a holomorphic extension $\G_\C\to M_\C$ of $\G\to M$, which is locally trivial with fiber type the complexification of the fiber type of $\G$.

\subsubsection{Holomorphic extension of compact Lie group actions}  
Let $G$ be a compact Lie group acting analytically on an analytic manifold $X$.  
Let $G_\C$ denote the complexification of $G$. By \cite{heinzner1994equivariantExt}, the action extends to a holomorphic action of $G_\C$ on a suitable complexification $X_\C$ of $X$.  
Consequently, the holomorphic action groupoid $G_\C \ltimes X_\C$ provides a holomorphic extension of the action groupoid $G \ltimes X$.

\section{Holomorphic extensions of s-proper groupoids}
\label{sec:s:proper}

Our aim in this section is to prove the following result.

\begin{theorem}
\label{thm:complex:s-proper:groupoid}
Every analytic s-proper Lie groupoid admits a holomorphic extension.
\end{theorem}

It is easy to give examples of analytic Lie algebroids which are not integrable. If an analytic algebroid $A\Rightarrow M$ is integrable, one may ask if the holomorphic extension $A_\C\Rightarrow M_\C$ is also integrable (as a germ around $M$). We don't know the answer to this question in general. 

\begin{remark}
     Since a holomorphic algebroid is integrable to a holomorphic groupoid if and only if the underlying real Lie algebroid is integrable (see \cite{laurent2009integration}), one could try to apply  the integrability criteria from \cite{crainic2003integrability} to answer this question. For that, one needs to compute the monodromy groups of the real algebroid underlying $A_\C\Rightarrow M_\C$ and show that they are uniformly discrete, using the assumption that the monodromy groupoids of $A\Rightarrow M$ are uniformly discrete. Recall from \cite{crainic2003integrability,CFM21} that the monodromy groups of $A$ are obtained from certain monodromy maps $\partial_x:\pi_2(\O,x)\to z(\mathfrak{g}_x(A))$, where $\O$ is the orbit of $A$ containing $x$. One then would hope to compare it with the monodromy map for the real algebroid underlying $A_\C$. The problem is that, as discussed in \cite{crainic2003integrability}, in general there is no way to control the monodromy maps when the point $x$ varies. Moreover, the real Lie algebroid underlying $A_\C$ can have  orbits contained in $M_\C - M$ but whose closure intersect $M$.  This may happen even if $A\Rightarrow M$ integrates to an s-proper groupoid (e.g., the cotangent Lie algebroid of the linear Poisson manifold $\mathfrak{so}^*(3)$). Hence, it seems hard to control the monodromy groups obstructing the integrability of $A_\C\Rightarrow M_\C$.
\end{remark}
 
 However, the previous theorem has the following corollary.

\begin{corollary}
If an analytic Lie algebroid admits an s-proper integration then it admits a holomorphic extension which is integrable.
\end{corollary}

\begin{proof}
    Let $\G\tto M$ be a Lie groupoid integrating an algebroid $A\Rightarrow M$. If $A$ is analytic then $\G$ is also analytic. This follows from the following two observations:
    \begin{enumerate}[(i)]
        \item It is enough to show that a neighborhood of the identity section $M\subset U\subset \G$ has a compatible analytic structure making $U$ a local analytic groupoid. Indeed, such a neighborhood generates $\G$, so one obtains a compatible analytic structure in the whole of $\G$. 
        \item The construction of the spray groupoid in \cite{CMS} shows that a choice of an analytic $A$-connection (which always exist) determines an analytic local groupoid $V$ integrating $A$. Possibly after restriction, $U$ and $V$ are isomorphic local Lie groupoids.
    \end{enumerate}
    Therefore, if $A\Rightarrow M$ has an s-proper integration $\G\tto M$, we can apply Theorem \ref{thm:complex:s-proper:groupoid} to obtain a holomorphic extension $\G_\C\tto M_\C$, whose Lie algebroid is a holomorphic extension of $A\Rightarrow M$.
\end{proof}

\begin{remark}
    It is also possible to show that analytic algebroids which (i) are transitive or (ii) have almost injective anchor, have integrable holomorphic extensions.
\end{remark}

\begin{corollary}
    Every analytic Poisson manifold of s-proper type admits an extension to an integrable holomorphic Poisson manifold.
\end{corollary}

\begin{proof}
    If $(M,\pi)$ is analytic of s-proper type then there exists an s-proper, analytic symplectic groupoid $(\G,\Omega)\tto M$ integrating it. In particular, the cotangent algebroid $(T^*M,\pi^\sharp,[\cdot,\cdot]_{\pi})$ has an s-proper integration. By the previous corollary, this algebroid admits a holomorphic extension which is integrable. Possibly after restriction, this holomorphic extension is the holomorphic cotangent algebroid of a holomorphic Poisson manifold $(M_\C,\pi_\C)$ extending $(M,\pi)$ -- see Example \ref{ex:Poisson:holomorphic}. 
\end{proof}

\begin{remark}
    In the previous proof, the symplectic groupoid $(\G,\Omega)\tto M$ has a groupoid holomorphic extension $\G_\C\tto M_\C$, but it is not clear it can be chosen to admit a holomorphic symplectic form $\Omega_\C$ extending $\Omega$. However, by the results of \cite{laurent2009integration}, the 1-connected integration of $T^*M_\C$ is a holomorphic symplectic groupoid integrating $(M_\C,\pi_\C)$.
\end{remark}

A similar result holds for analytic Dirac structures of s-proper type.

\subsection{The strategy to prove Theorem  \ref{thm:complex:s-proper:groupoid}}
The proof Theorem \ref{thm:complex:s-proper:groupoid} consists of the following steps:
\begin{enumerate}[Step 1.]
\item Analytically linearize $\G\tto M$ on invariant opens $U_i$ around each of its orbits;
\item Construct full holomorphic extensions $(\G_i)_\C\tto (U_i)_\C$ for the restrictions $\G_i:=\G|_{U_i}$ using the linear local model;
\item Assemble the full local holomorphic extensions $(\G_i)_\C$ into a holomorphic extension of $\G$.
\end{enumerate}
Note that Step 1 follows immediately from Theorem \ref{thm:analytic_linearization}. The next two paragraphs implement the other two steps.

\subsection{Holomorphic extension of the linear local model}

Assume that $\G\tto M$ is an analytic s-proper groupoid and let $\O$ be an orbit of $\G$. Fix $x\in \O$ and let $P:=s^{-1}(x)$, so $t: P \rightarrow \O$ is a principal bundle over $\O$ with structure group $G:=s^{-1}(x)\cap t^{-1}(x)$. The group $G$ acts on the normal space  $N=T_x M/T_x\O$ and the linear local model for $\G$ can be recast as follows (see, e.g., \cite{crainic2013linearization}):
\begin{itemize}
    \item The normal bundle of $\O$ is isomorphic to the associated vector bundle 
    \[ \N_\O \cong (P\times N)/ G,\]
    where the $G$-action on $P \times N$ is $g (p,v) = (pg^{-1}, g v)$.
    \item The normal bundle of $\G_\O$ is isomorphic to
     \[  \N_{\G_\O}\tto \N_\O  \cong (P\times P\times N)/ G\tto (P\times N)/G. \]
    where the $G$-action on $P\times P\times N$ is $g(p_1,p_2,v):=(p_1g^{-1},p_2g^{-1},gv)$.
    \end{itemize}
Viewing $P\times P\times N\tto P\times N$ as the groupoid product of the pair groupoid $P\times P\tto P$ and the identity groupoid $N\tto N$, the $G$-action is by groupoid automorphisms and the resulting groupoid structure on the quotient is isomorphic to the linear local model:
\begin{equation}
\label{eq:local:normal:form}
\vcenter{
\xymatrix@R=5pt{\N_{\G_\O} \ar@<0.25pc>[dd] \ar@<-0.25pc>[dd]  & & (P\times P\times N)/ G \ar@<0.25pc>[dd] \ar@<-0.25pc>[dd] \\ & \cong & \\ \N_\O & & (P\times N)/G}}
\end{equation}

Letting  $N_\C:=N\otimes\C$ and denoting by $P_\C\to \O_\C$ the principal $G_\C$-bundle extending $P\to \O$, given by Theorem \ref{thm:complexify_principal_bundle}, we now have the following result.

\begin{proposition}
The linear local model around an orbit \eqref{eq:local:normal:form}  admits a full holomorphic extension of the form
 \begin{equation}
 \label{eq:hol:extension:linear model}
 (P_\C\times P_\C\times N_\C)/G_\C\tto (P_\C\times N_\C)/G_\C.
 \end{equation}
\end{proposition}

\begin{proof}
The proof is a generalization of Example \ref{sec:complexify:proper groupoid}. First, we have the obvious embedding 
\[ \iota:(P\times P\times N)/ G\tto (P\times N)/ G\hookrightarrow (P_\C\times P_\C\times N_\C)/G_\C\tto (P_\C\times N_\C)/G_\C, \] 
and property (C1) clearly holds. 

To check the other two properties, observe that for any holomorphic groupoid $\H\tto X$, any morphism 
\[ \phi:(P\times P\times N)/ G\to \H \] 
can be written in the form
\[ \phi([p_1,p_2,v])=\psi(p_1,v)\psi(p_2,v)^{-1}, \]
for a map $\psi:P\times N\to s^{-1}(x)$ satisfying:
\begin{enumerate}[(a)]
\item $\psi$ is $G$-equivariant relative to a Lie group morphism $\rho:G\to \H_x$;
\item $\psi$ covers a map $\psi_0:(P\times N)/G\to X$.
\end{enumerate}

If $\phi$ is analytic, so is $\psi$. In this case, it follows from \cite{heinzner1994equivariantExt} that $\psi$ has a holomorphic extension $\psi^*:Q\to s^{-1}(x)$ defined in some $G_\C$-invariant neighborhood $Q$ of $P\times N$ in $P_\C\times N_\C$. This extension is $G_\C$-equivariant relative to the complexification of the morphism $\rho:G\to \H_x$. 

Next, we have the $G_\C$-equivariant embedding
\[ Q\times_{N_\C} Q\to P_\C\times P_\C\times N_\C, \quad (p_1,v,p_2,v)\mapsto (p_1,p_2,v), \]
and we can define a holomorphic groupoid morphism
\[ \phi^*:Q\times_{N_\C} Q/G_\C\to \H, \quad  \phi^*([p_1,p_2]):=\psi^*(p_1)\psi^*(p_2)^{-1}. \]
Notice that $Q\times_{N_\C} Q/G_\C \tto Q/G_\C$ is a full groupoid neighborhood of $(P\times P\times N)/ G$ in $(P_\C\times P_\C\times N_\C)/G_\C$, and properties (C2) and (C3) now follow.
\end{proof}

\subsection{Gluing the local holomorphic extensions}

Step 3 of the proof will follow from the previous two steps and the following general result.

\begin{proposition}
Let $\G\tto M$ be an analytic Lie groupoid and assume that there is an open $\G$-invariant cover $\{U_i\}_{i\in I}$ of $M$ such that each restriction $\G|_{U_i}\tto U_i$ admits a full holomorphic extension. Then $\G\tto M$ admits a holomorphic extension.
\end{proposition}

\begin{proof}
Note that the orbits of an s-proper Lie groupoid are compact. Hence, we can assume that the open $\G$-invariant cover $\{U_i\}_{i\in I}$ is locally finite and that the closures $\overline{U_i}$ are compact. Let $M_\C$ be a complexification of $M$, and let $\G_i^*\tto U^*_i$ be a full holomorphic extension of $\G|_{U_i}\tto U_i$. Then both $M_\C$ and $U^*_i$ are complexifications of $U_i$, hence by restricting $\G_i^*\tto U^*_i$ we can assume that $\{U^*_i\}_{i\in I}$ is a locally finite collection of open sets in $M_\C$ satisfying
\[ U_i=U^*_i\cap M,\quad  \overline{U_i}=\overline{U^*_i}\cap M. \] 

Set $U_{ij}:=U_i\cap U_j$ and $U^*_{ij}:=U^*_i\cap U^*_j$. The restrictions $\G_i^*|_{U^*_{ij}}$ and $\G_j^*|_{U^*_{ij}}$ are both full holomorphic extensions of $\G|_{U_{ij}}$. It follows that one can choose $T^*_{ij}$, open neighborhoods of $U_i\cap U_j$ in $M_\C$, satisfying 
\[ \overline{U_i}\cap\overline{U_j}\subset T^*_{ij}, \]
for which there exist holomorphic groupoid isomorphisms
\[ \psi^*_{ij}:\G_i^*|_{T^*_{ij}}\diffto \G_j^*|_{T^*_{ij}}, \quad \psi^*_{ji}=(\psi^*_{ij})^{-1}. \]

Applying Lemma \ref{lemma:refine:cover}, one can find a family of open sets $\{W^*_i\}_{i\in I}$ in $M_\C$ such that:
\begin{enumerate}[(i)]
\item $W^*_i\subset U_i^*$ 
\item $W^*_i\cap M=U_i$;
\item If $U_i^*\cap U^*_j \neq \emptyset$, then $W^*_i\cap W^*_j \subset T^*_{ij}$.
\end{enumerate}
Then the restrictions $\G_i|_{W^*_i}$ form a collection of full holomorphic extensions of the $\G_i$'s and the $\psi^*_{ij}$'s restrict to holomorphic groupoid isomorphisms
\[ \psi^*_{ij}:\G_i^*|_{W^*_i\cap W^*_j}\diffto \G_j^*|_{W^*_i\cap W^*_j}, \quad \psi^*_{ji}=(\psi^*_{ij})^{-1}. \]
Replacing $M_\C$ by the union of the $W^*_i$'s, we obtain a holomorphic groupoid 
\[ \G_\C:=\frac{\bigsqcup_i \G^*_i|_{W^*_i}}{\sim}\tto M_\C,\]
where $g\sim \psi_{ij}(g)$. 

It remains to prove that $\G_\C\tto M_\C$ satisfies properties (C1)-(C3) in Definition \ref{def:holomorphic:extension:groupoid}. The obvious embedding $\iota:\G\to\G_\C$ is totally real so (C1) is obvious. Let $\phi:\G\to\H$ be real analytic morphism into some holomorphic groupoid $\H$. The restriction of $\phi$ to each $\G_i$ has a holomorphic extension $\phi^*_i$ defined on some full groupoid neighborhood $\G_\C|_{V^*_i}$ of $\G_i$. Eventually after restriction of the $V^*_i$'s, one can assume that both the $V^*_i$'s and the intersections $V^*_i\cap V^*_j$ are connected. Since the $\G_i$'s have connected s-fibers, it follows that $\phi^*_i=\phi^*_j$, by uniqueness of holomorphic extension. Hence, letting $V^*:=\bigcup_{i\in I} V^*_i$, we can define a holomorphic groupoid morphism 
\[ \phi^*:\G_\C|_{V^*}\to \H, \quad \phi^*|_{\G_{V^*_i}}:=\phi^*_i, \]
which extends $\phi$. This shows that (C2) holds. Condition (C3) follows from the fact that (C3) is satisfied for the local holomorphic extensions $\G_i^*$.
\end{proof}

\appendix

\section{Analytic Haar densities}
\label{sec:analytic_haar_density}
In this section, we will show that any analytic s-proper groupoid $\G$ admits an analytic Haar density.

Let $E$ be a vector bundle on $M$. We denote by 
\[ \D(E):= \bigsqcup_{p \in M} \D(E_{p}) \]
the density line bundle whose sections are the densities of $E$. In this paper we only consider \emph{positive densities}, i.e., densities $\mu$ such that  $\mu(e_1,\cdots,e_n)>0$ for any  local frame $\{e_1,\cdots,e_n\}$ for $E$. When the bundle $E$ is analytic, the line bundle $\D(E)$ is analytic and we say that a density $\mu$ is \textbf{analytic} if it is an analytic section of $\D(E)$.

For a manifold $M$ we let  $\D(M):= \D(T^{*}M)$ and we call sections of $\D(M)$ densities on $M$. When the manifold $M$ is oriented, a volume form $\omega$ induces a positive density $\mu$ by setting $\mu_p:= |\omega_p|$ for every $p$ in $M$. Not every manifold is oriented and admits a volume form, but every manifold admits a density and we can further assume it to be analytic when $M$ is analytic. In fact, we have:


\begin{proposition}
    Every analytic vector bundle has a positive analytic density.
\end{proposition}

To prove this, we recall the following well-known criterion for analytic functions \cite{krantz2002primer}.
\begin{lemma}\label{analyticity_criterion}
A smooth function $f: U \rightarrow \R$ is an analytic function on an open set $U\subset \R^n$ if and only if for any compact set $K \subset U$, there exists a constant $C>0$ depending on $K$ such that for all $\alpha\in \mathbb{Z}_{\geq 0}$ and all $x\in K$,
\[ \left| \frac{\partial^\alpha f}{\partial x^\alpha} \right| \leq C^{\alpha+1} \alpha !\]
\end{lemma}

By the classical results of Grauert \cite{Grauert58} and Morrey \cite{Morrey58} every analytic bundle $E$ admits an analytic Riemannian metric $g$. It has an associated positive density $\mu$ which on a local frame $\{e_1,\cdots,e_n\}$ for $E$, where $g=g_{\alpha \beta} e_\alpha\otimes e_\beta$, is given by
\[ \mu = \sqrt{\det(g_{\alpha \beta})}|e_1 \wedge \cdots \wedge e_n|. \] 
This density is analytic because, by the above criterion for analyticity, $\sqrt{f}$ is analytic whenever $f$ is a positive analytic function.


\begin{definition}\label{def:analytic haar density}
Let $\G\tto M$ be an analytic s-proper groupoid. We say that a family of positive densities $\{\mu^x\}_{x\in M}$ on the $s$-fibers of $\G$ is an \textbf{analytic normalized Haar density} if the following conditions hold:
\begin{enumerate}[(i)]
    \item \emph{Analyticity}: The function 
    \[x \mapsto \int_{s^{-1}(x)} f(g)\mu^x(g)\]
is analytic on $M$ for any analytic function $f\in C^{\omega}(\G)$.
\item  \emph{Right-invariance}: For any arrow $y \xleftarrow{h} x$ and $f \in C^\omega(s^{-1}(x))$, we have
    \[
    \int_{s^{-1}(y)} f(gh) \mu^y(g) = \int_{s^{-1}(x)} f(g) \mu^x(g).
    \]
\item \emph{Normalization}: For all $x\in M$
\[\int_{s^{-1}(x)} \mu^{x}(g) = 1.\]
\end{enumerate}
\end{definition}

The rest of this section is devoted to the proof of the following fact.

\begin{theorem}\label{thm:analytic haar density}
    Every $s$-proper groupoid $\G$ admits an analytic normalized Haar density.
\end{theorem}

Let $A$ be the Lie algebroid of $\G$ and let $\mu$ be any positive analytic density on the dual bundle $A^*$. Consider the 
vector bundle map
\[
\vcenter{\xymatrix{
\Ker(\d s) \ar[r]^{\phi} \ar[d] & A \ar[d] \\
\G \ar[r]_{t} & M
}}\quad \phi(v):= \d R_{g^{-1}}(v). 
\]
Pulling back $\mu$ through $\phi$ yields an analytic density $\Tilde{\mu}$ on $(\Ker(\d s))^*$, i.e., a family of analytic densities $\{\Tilde{\mu}^{x}\}$ on each fiber $s^{-1}(x)$. It is easy to check that for any arrow $y \xleftarrow{h} x$, $\Tilde{\mu}^{y} =R^{*}_h \Tilde{\mu}^{x}$, so $\{\Tilde{\mu}^x\}$ is right-invariant. By Lemma \ref{analyticity}, the positive function $c:\G\to\R$ defined by
\[c(x):= \int_{s^{-1}(x)} \Tilde{\mu}^{x}(g).\]
is analytic. Replacing $\Tilde{\mu}^{x}$ by $\frac{1}{c(x)} \Tilde{\mu}^{x}$, we obtain an analytic normalized density. Finally,  property (i) in Definition \ref{def:analytic haar density} follows from the following useful lemma.

\begin{lemma}\label{analyticity}
Let $s: E\rightarrow M$ be a proper analytic submersion. Given an analytic density $\mu$ on the vector bundle $(\Ker(\d s))^* \rightarrow E$, for any analytic function $f:E\to\R$ the average
 \[x \mapsto \int_{s^{-1}(x)} f(g)\mu^x(g)\]
is an analytic function on $M$. 
\end{lemma}

\begin{proof}
By the analytic version of Ehresmann theorem, $s$ is a locally trivial fibration. Choosing some open neighborhood $U$ of $x$, we can assume that $E|_{U} = U \times s^{-1}(x)$ and $s = \text{pr}_{U}$. Let $\{V_i\}$ be a finite open cover of $s^{-1}(x)$ and $\{\rho_i\}$ be a partition of unity subordinate to it.

\[\int_{s^{-1}(x)} f(g)\mu^x(g) =\sum_i  \int_{ V_i} \rho_i f(g) \mu^x(g) .\] 
It suffices to prove that $f_i(x): = \int_{s^{-1}(x)\cap V_i} \rho_i f(g) \mu^x(g)$ is analytic. In local coordinates, write $g = (x,y)$ and $\mu^{x}(g) = \lambda(x,y)dy$, we have
\[f_i(x) = \int_{V_i} \rho_i(y)f(x,y)\lambda(x,y)dy. \]

Since $f$ and $\lambda$ are analytic, so is $f\lambda$. By Lemma \ref{analyticity_criterion}, for any compact set $K \subset U \times V_i$, there exists constant $C$ such that $\left|\frac{\partial^\alpha (f \lambda)}{\partial x^\alpha}(x,y)\right| \leq C^{|\alpha|+ 1}\alpha!$. Therefore
 \[ \left|\frac{\partial^\alpha f_i}{\partial x^\alpha} (x) \right|\leq  \int_{\text{supp}(\rho_i)} \left|\frac{\partial^\alpha (f \lambda)}{\partial x^\alpha}(x,y)\right|dy  \leq |\text{supp}(\rho_i)|C^{|\alpha|+ 1}\alpha !,\]
 proving that $f_i$ is analytic.
\end{proof}

\section{Complexification and refinement of covers}
\label{sec:covers}

Let $M$ be an analytic manifold and $M_\C$ a complexification of $M$. We identify $M$ with its image under the closed embedding $\iota:M\to M_\C$. We need the following lemma concerning refining complexifications of covers. Similar results were used in the proof of existence of holomorphic extensions by Whitney and Bruhat \cite{whitney1959quelques}).

\begin{lemma}
\label{lemma:refine:cover}
Let $\{U_i\}_{i\in I}$ be a locally finite cover of $M$ by open sets with compact closures and let $\{U^*_i\}_{i\in I}$ be a locally finite family of open sets of $M_\C$ such that:
\[ U_i=U^*_i\cap M,\quad  \overline{U_i}=\overline{U^*_i}\cap M. \] 
Suppose one is given open neighborhoods $T^*_{ij}$ of $U_i\cap U_j$ in $M_\C$ satisfying
\[ \overline{U_i}\cap\overline{U_j}\subset T^*_{ij}. \]
Then one can find a family of open sets $\{W^*_i\}_{i\in I}$ in $M_\C$ such that:
\begin{enumerate}[(i)]
\item $W^*_i\subset U_i^*$ 
\item $W^*_i\cap M=U_i$;
\item If $U_i^*\cap U^*_j \neq \emptyset$, then $W^*_i\cap W^*_j \subset T^*_{ij}$.
\end{enumerate}
\end{lemma}

\begin{proof}
For any $i,j\in I$, let 
   \[ V^*_{ij} := (U^*_i-U^*_i \cap U^*_j) \cup T^*_{ij}.\]
We claim that the interior  $\mathring{V}^*_{ij}$ is an open neighborhood of $U_i$. To see this, given $x\in U_i$, we consider two possibilities:
\begin{enumerate}[(a)]
\item $x\in \overline{U_i}\cap\overline{U_j}$: then $x\in T^*_{ij}$, so it is an interior point of $V^*_{ij}$. 
\item $x\not\in  \overline{U_i}\cap\overline{U_j}$: since $\overline{U^*_i}\cap M\subset\overline{U_i}$, we have
\[ \overline{U^*_i\cap U^*_j} \cap M\subset\overline{U^*_i}\cap\overline{U^*_j} \cap M\subset\overline{U_i}\cap\overline{U_j} \]
Hence, $x$ is an interior point of $U^*_i-U^*_i \cap U^*_j$ and, hence, also of $V^*_{ij}$.
\end{enumerate}
We now define
 \[W^*_i := \left(\bigcap_{U^*_i\cap U^*_j\neq \emptyset} \mathring{V}^*_{ij} \right)\cap U^*_i. \]
Clearly, $W^*_i$ is an open set contained in $U^*_i$ which, by the previous claim, contains $U_i$. Hence, $W^*_i\cap M=U_i$. Finally,  If $U^*_i\cap U^*_j \neq \emptyset$, then
 \[ W^*_j\cap W^*_j \subset \mathring{V}^*_{ij}\cap U^*_i\cap U^*_j \subset T^*_{ij}.\]
\end{proof}

As an application of the previous lemma, we give a proof of Theorem \ref{thm:complexify_principal_bundle} different from the original one in \cite{heinzner1994equivariantExt}. 

Assume that $P$ is an analytic principal $G$-bundle over an analytic manifold $M$, with $G$ a compact Lie group. Choose a trivializing open cover $\{U_i\}_{i\in I}$ of $M$, which we can assume is locally finite, has compact closures and connected double intersections. Denote by $g_{ij}:U_i\cap U_j\to G$ the corresponding transition functions. 

Let $M_\C$ be a complexification of $M$ and choose $\{U^*_i\}_{i\in I}$ a locally finite family of open sets of $M_\C$ extending the family $\{U_i\}_{i\in I}$, as in Lemma \ref{lemma:refine:cover}, with connected double intersections (see Remark \ref{rem:extending:cover}). We can holomorphically extend the transition functions to
\[ g^*_{ij}: T^*_{ij} \rightarrow G_\C\]
where the $T^*_{ij}$ are open neighborhoods of $U_i\cap U_j$ in $M_\C$. By analytic continuation, possibly after shrinking the $T^*_{ij}$, we may also assume that the family $\{g^*_{ij}\}$ satisfies the cocycle condition. 

Now let $\{W^*_i\}$ be the family of open sets given by Lemma \ref{lemma:refine:cover}. We can restrict the transition function $\tilde{g}_{ij}$ to $W^*_i\cap W^*_j$ and define:
\[ P_\C = \bigsqcup_i W^*_i\times G_\C/\sim\]
where $(x,\lambda)\sim (x,g_{ij}\lambda)$. Replacing $M_\C$ by the union of the $W^*_i$'s, we obtain a holomorphic principal $G_\C$-bundle $P_\C\to M_\C$ which is easily seen to satisfy the properties in Theorem \ref{thm:complexify_principal_bundle}.

\end{document}